\documentclass[12pt, a4paper,twoside,openright]{article}

\usepackage{latexsym}
\usepackage{indentfirst}
\usepackage{graphicx}
\usepackage{amsmath,amssymb,amsfonts,amsthm,amstext,amscd}
\usepackage{mathrsfs}
\usepackage{fancyhdr}
\usepackage{times}
\usepackage{enumerate}
\usepackage{natbib}
\usepackage[hang,flushmargin]{footmisc}

\fancypagestyle{plain}{%
\fancyhf{}%
\fancyhead[LO]{\textbf{Science \& Philosophy}
}
\fancyhead[RO]{Volume 12(1), 2024 }

\fancyfoot[CO,RE]{}
}

\begin{document}
\title{\textbf{Logic as an internal organisation of language}}
\author{Boris Čulina\footnote{University of Applied Sciences Velika Gorica, Velika Gorica, Croatia; boris.culina@vvg.hr}}

\date{}

\maketitle

\begin{abstract}
\begin{normalsize}
\noindent
Contemporary semantic description of logic is based on the ontology of all possible interpretations, an insufficiently clear metaphysical concept. In this article, logic is described as the internal organization of language.  Logical concepts -- logical constants, logical truths and logical consequence -- are defined using the internal syntactic and semantic structure of language. For a first-order language, it has been shown that its logical constants are connectives and a certain type of quantifiers for which the universal and existential quantifiers form a functionally complete set of quantifiers. Neither equality nor cardinal quantifiers belong to the logical constants of a first-order language. \footnote{\noindent Received on March 14th, 2024. Accepted on Jule 6th, 2024. Published on xxx, 2024. doi: xxx. ISSN 2282-7757; eISSN 2282-7765. \copyright The Authors. This paper is published under the CC-BY licence agreement.}
\end{normalsize}
\end{abstract}

%
%

\newpage

\section{Introduction}

Starting with Tarski's definition of the concept of logical consequence from 1936 (\cite{Tar1}) and the definition of the concept of logical constant from 1966 (\cite{Tar2}), the examination of logical concepts is dominated by an approach in which logical concepts are concepts that are in some way invariant to language interpretations.\footnote{The basics of this approach can already be found in Bolzano from 1837 (\cite{Bol}).} There are differences in how the concept of invariance is formulated. See, e.g., \cite{Bon} for an overview and defence of such an approach in the case of the concept of logical constant. Of course, there are other approaches,\footnote{See, e.g., \cite{Mac} for an overview in the case of the concept of  logical constant.} but this approach is dominant and here I will call it \textit{the received view}.  The received view, no matter how it is formulated, is based on the ontology of all possible interpretations, an insufficiently clear metaphysical concept. In mathematical logic, interpretations are usually realized in the world of sets: sets, relations, and functions are the values of the interpretations. However, the set theory based on ZFC axioms does not have a clear intuitive basis to give us a satisfactory answer to the question of which sets exist and which do not exist.\footnote{See, e.g., \cite{Cul}.} Consequently, interpretations in set theory are also insufficiently clear. Ultimately, these interpretations are translations of the language whose interpretations we examine into the language of sets. In this translation, for example, the logical constants of the translated language are described by the same logical constants of the language of sets, which does not contribute to their better understanding. Thus, in my opinion, the received view, due to the use of vague metaphysical assumptions, is unacceptable.

This article is based on an understanding of logical concepts -- logical constants, logical truths, and logical consequence -- that is  opposite  to the received view. If we were to call the understanding of logical concepts in the received view an external understanding then this would be an internal understanding of logical concepts. The analysis of logical concepts that will be conducted here is based on the assumption that logic is always the logic of a language, how we apply the language, and how its parts are syntactically and semantically connected. It means that logic is not an objective science about propositions, thoughts, absolute truths or some other universal metaphysical ghosts, but it is a normative science of the organization and use of a particular language. The logic of a language is just the inner organization of the language together with external assumptions of its use. In this article I will not defend this understanding of logic. This was done in article \cite{Cu0}, where the essential role of language in our rational cognition and thinking was analysed. In this article, I will use the example of a first-order language to demonstrate how logical concepts can be described using this understanding of logic.

\section{What is logical in first-order logic?}

The external assumptions of an interpreted first-order language are: (i) the language has its own domain of interpretation -- a collection of objects that the language speaks of, (ii) every constant denotes an object, and every variable in a given valuation denotes an object, (iii) every function symbol symbolizes a function which applied to objects gives an object,  (iv) every predicate symbol symbolizes a predicate which applied to  objects gives a truth value, $ True $ or $ False $.  The fulfilment of the external assumptions is crucial for the application of a language but not for the logic of the language. The only important thing for the logic of the language is that these assumptions are part of the specification of the language, not whether they are fulfilled. Regardless of whether the external assumptions are fulfilled or not, the logic of the language demands that when we use the language we assume that they are fulfilled. In thinking itself there is no difference whether we think of objects that really exist or we think of objects that do not really exist and whether the predicate symbols we use can be applied to such objects at all or not. That difference can be registered only in a ``meeting'' with reality. 
Thus, external assumption of the language use have no ontological weight for the logic of the language.  

The inner organization of a first-order language   is determined by the rules of the construction  of more complex language forms from simpler ones, starting with names, variables and function symbols for building terms, and with atomic sentences for building sentences. In these constructions we use special symbols which identify the type of the construction. With each construction, and thus the symbol of the construction, a semantic rule is associated that determines the semantic value of the constructed whole using the semantic values of the parts of the construction.\footnote{In a given interpretation and a given valuation of variables, the semantic value of a term is the object described by the term and the semantic value of a sentence is its truth value.} The symbol of a language construction will be termed \textit{logical symbol} or \textit{logical constant} if the associated semantic rule  is an internal language rule: the rule does not refer to the reality the language speaks of, except possibly referring to external assumptions of the language use. It may be objected that this description of the concept of  logical symbol is not clear enough. But an incomparably greater ambiguity is present in the received view where the description of logical constants encompasses all interpretations of a language. While in the received view further refinement of the term ``all interpretations'' necessarily involves metaphysical assumptions, it will be shown below that the description proposed here is self-sufficient and clear enough to give us the answer in a concrete situation whether the symbol of a construction is a logical symbol or not. Since this is an internal language approach in determining logical symbols -- an approach that, in addition to external assumptions about the language use, does not include the reality of which the language speaks -- it is automatically topic-neutral. Everything defined in it as a logical concept will be a logical concept in  the received view approach. Thus the received view gives only the necessary conditions for logical symbols. 

The description of logical concepts in terms of the internal structure of language is present in the literature in various forms. However, as far as I know, it has neither been given sufficient importance nor has it been pronounced precisely enough. For example, in the review article \cite{Sha} on logical consequence, giving various criteria for the notion of logical consequence, Shapiro also mentions the following criterion:

\begin{quote}
	
	$\Phi$ is a logical consequence of $\Gamma$ if the truth of the members of $\Gamma$ guarantees the truth of $\Phi$ in virtue of the meanings of the logical terminology.
\end{quote}	

\noindent However, what ``the meaning of logical terminology'' means is not specified. Likewise, Quine in \cite{Qui}, page 48, among other criteria, gives the following criterion for logical consequence:

\begin{quote}
	One closed sentence logically implies another when, on the assumption that the one is true, the structures of the two sentences assure that the other is true. The crucial restriction here is that no supporting supplementary assumption or information be invoked as to the truth of additional sentences. Logical implication rests wholly on how the truth functions, quantifiers, and variables stack up. It rests wholly on what we may call, in a word, the logical structure of the two sentences. 
\end{quote}	

\noindent Although Quine is known for precision, this description is also not precise enough nor is it further specified. That Quine's approach is different from the approach in this article can also be seen below from the consequences derived. In contrast to the above and other similar descriptions known to me, in this paper the semantic extensional rules of the first-order language constructions are precisely specified. They give a clear criterion whether a construction symbol is a logical symbol or not, whether a sentence is a logical truth and whether a sentence follows logically from a set of sentences. An analysis of logical symbols of a first-order language follows.

Connectives of a first-order language give us one  way to combine simpler sentences into more complex ones. Every connective, regardless of whether it is abstracted from a corresponding natural connective or not, is determined by a Boolean function $f: \{True, False\}^n\ \longrightarrow \ \{True , False\}$, where $ n $ is non-negative integer. This function  describes how the truth value of the sentence composed by this connective depends on the truth values of sentences from which it is composed. Since Boolean functions are internal semantic functions of the language, functions independent of the reality the language speaks of, these connectives are logical symbols of the language. Of course, the well-known results on functionally complete sets of logical connectives show that in a first-order language we should not have connectives other than standard ones, for example, $\land$, $\lor$, $\lnot$, $\rightarrow$ and $\leftrightarrow$. All other connectives can be defined using these.

Qualitatively different  way of combining sentences to more complex sentences is by combining with quantifiers ``for all'' and ``exists'', which are symbolized by symbols $ \forall $ and $ \exists $. How to characterize this type of combination? Are there other quantifiers of this type? Can we express all of them by quantifiers $ \forall $ and $ \exists $ which are abstracted from natural language?  The approach conducted here is inspired by the description of the universal quantifier in \cite{GB}, page 40.

Let's take, for example, the quantifier $ \forall $. From a sentence $S(v)$, where $v$ is a variable, not necessarily free in the sentence, by the symbol $ \forall $ the more complex sentence $ \forall v \ S(v) $ is built. In a given interpretation of the language  and in a given valuation of all variables except $v$, the truth value of the sentence $S(v)$ depends on the valuation of the variable $v$. To determine the truth value of $ \forall v \ S(v) $, we must determine the truth values of  $S(v)$ for all valuations of $v$. We can get three sets of  truth values: $ \{ True\} $ (for all valuations of $v$ the sentence $S(v)$ is true), $ \{ False\} $ (for all valuations of $v$ the sentence $S(v)$ is false) and $ \{ True,False\} $ (for some valuations of $v$ the sentence $S(v)$ is true and for some valuations it is false). In the first case the sentence $ \forall v \ S(v) $ is true, in other cases it is false. So, the quantifier $ \forall $ is determined by the function that maps non-empty sets of truth values to truth values. The quantifier $ \exists $ is of the same type: it is determined by the function that maps all  sets which contain $ True $  to $ True $, and $ \{ False\} $ maps to $ False $. Every such function that maps non-empty set of truth values to truth values determines an extensional construction. These functions will be termed \textit{logical quantifier functions}. The corresponding syntactical symbol of the construction will be termed \textit{logical quantifier}. Just as logical connectives are logical symbols because they are determined by internal functions that map truth values to truth values, so logical quantifiers are logical symbols because they are determined by internal functions that map sets of truth values to truth values. It will be shown below that these are the only logical quantifiers of a first-order language, which will justify the name given to them: ``logical quantifier''.\footnote{In this paper, only type $<1>$ quantifiers will be analysed (see, e.g., \cite{We}).}

Since there are $ 2 ^ 3 = 8$ functions from the set of non-empty sets of truth values to the set of truth values, there are 8 logical quantifiers. However, they do not need to be introduced by special constructs in a first-order language because all the others can be defined using $\forall $ and $\exists $. Definitions are given in the following table ($\top$ is the label for \textit{True}, $\bot$ is the label for \textit{False}):

\begin{center}
	\includegraphics[height=5cm]{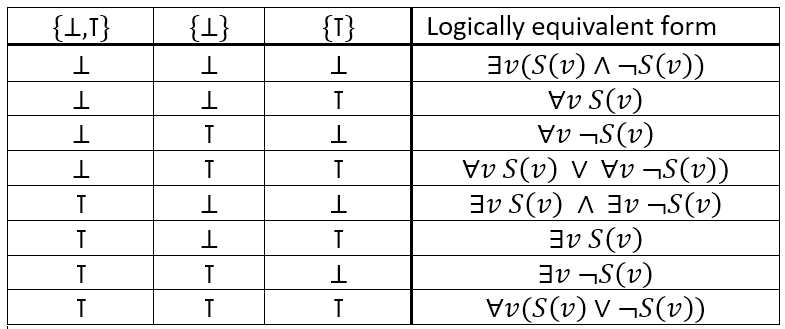}
	
	Figure 1: Logical quantifiers
\end{center}

\noindent It means that the set $\{\forall, \exists \}$ is a functionally complete set of logical quantifiers of a first-order language. Because $ \forall v \ S(v) $ is logically equivalent to  $ \lnot\exists v\lnot \ S(v) $ and $ \exists v \ S(v) $ is logically equivalent to  $ \lnot\forall v\lnot \ S(v) $, sets $\{\forall \}$ and $\{\exists \}$ are also  functionally complete sets of logical quantifiers of a first-order language. Analogous to the results for connectives, this result shows that in a first-order language we do not need to have other logical quantifiers besides the standard ones, $\forall$ and $\exists$.

Is equality (the symbol $ = $) a logical symbol of a first-order language? Let $ t_1 $ and $ t_2 $ be terms. Using the symbol $ = $ the atomic sentence $ t_1 = t_2 $ is built. The corresponding semantic construction is as follows: in a given interpretation and a given valuation of variables, the sentence $ t_1 = t_2 $ is true iff $ t_1 $ and $ t_2 $ denote the same object. This is a construction that maps the semantic values of these two terms -- denoted objects -- to the truth value of the corresponding atomic sentence. However, this construction delves deeper into reality than the external assumption of the language which only says that each term  denotes an object. To determine whether terms denote the \textit{same} object, we must look at the reality that language speaks of. For example, the key to the  Superman story is the claim that Superman = Clark Kent. To determine whether this is true or not, logic does not help us but we have to look into the comics reality. Or Frege's example: ``The morning star = The evening star''. We know which objects are denoted by these terms, but logic is not enough to determine that it is the same object -- we need astronomical observations. Determination of the truth value of the sentence $ t_1 = t_2 $ using the semantic values of the terms $ t_1 $ and $ t_2 $ involves reality beyond external assumptions of the language use. Hence, equality is not a logical symbol. The reason why the logicalness of equality is the subject of dispute\footnote{Quine’s doubts and pro et contra arguments about the logical status of equality can be seen in \cite{Qui}, pages 61--64.}  may lie in the fact that, unlike, for example, the comparison of numbers, we can state certain logical truths about equality. This is because the description of the symbol of equality includes the internal structure of language  in one part -- it mentions the denotations of terms. For example, the external assumption of the language use is that in a given valuation, each variable  denotes a specific object. So $ x $ and $ x $ will denote the same object. Therefore  $ \forall x \ x = x $ is a logical truth. 

Since equality is not a logical symbol, quantifiers that are necessarily described by equality are also not logical symbols. Such is, for example, the quantifier ``there is one and only one'' (the symbol $ \exists !$):

$$\exists !v  \ S(v) \leftrightarrow \exists x S(x) \land \forall y\forall z(S(y) \land S(z) \rightarrow y=z)$$ 

\noindent We can also see in a direct way that such a quantifier is not a logical symbol -- by examining the semantic rule of the associated language construction $S(v) \mapsto \exists !v  \ S(v)$. In a given interpretation and valuation of all variables except $ v $, the sentence  $\exists !v  \ S(v)$ is true iff the sentence  $S(v)$ is true in exactly one valuation of $ v $. We can describe this construction by a function that maps multisets composed of truth values to truth values.\footnote{A multiset composed of truth values is a function from the set $\{True, False\}$ into the set of non-negative integers. The value of a multiset on a given truth value is called its multiplicity.} If the multiplicity of $ True $ is equal to 1, the function gives the value $ True $, otherwise it gives the value $ False $. This function, like the Boolean and quantifier functions, is an internal semantic function, a function that connects semantic values independently of the reality the language speaks of. However, the overall semantic rule of this construction includes the reality because the argument of the function, a multiset, cannot be formed without distinguishing objects from the reality. How many times a truth value has occurred cannot be determined without distinguishing valuations of $ v $, that is, without determining when a valuation yields the same object and when it does not. And this requires, as with equality, knowledge of the reality the language speaks of, knowledge which goes beyond the external assumptions of the language use. This argument is easy to generalize. All cardinal quantifiers  are not logical quantifiers, because the semantics of the language construction determined by such a quantifier is described by the same type of function as for $ \exists! $ -- by a function that maps multisets composed of truth values to truth values. With all these quantifiers, in a given interpretation and valuation of variables, the identification of the multiset on which the function acts includes reality in the same way as with $ \exists! $. So these quantifiers are not logical symbols.

According to the standard recieved view that is clearly presented in \cite {Sher}, logical constants are constants that are invariant under bijections between domains. They include all  cardinal quantifiers, including infinite cardinals, which state how many objects satisfy a formula. In response to such a broad concept of logical constant, in which the distinction between mathematics and logic is lost, the articles \cite {Fef} and \cite {Bon1} have emerged that widen the conditions of invariance, thus narrowing the concept of logical constant. In the criterion of invariance, Feferman replaces the notion of bijection with the notion of surjection, and this substitution leads precisely to the logical constants of  first-order languages established in this article.

Like the concept of  logical symbol, so the concepts of logical truth and logical consequence have a basis in the internal structure of a language. \textit{Logical truth} is a sentence of a language that is true not in terms of what reality is but in terms of a language we use to describe reality: it is truth determined by the internal semantic structure of the language. For example, $ \lnot A\ \lor A $ is a logical truth of the first-order languages, because its truth is determined by the internal structure of the language, in this case the semantics of the connectives $ \lnot $ and $ \lor $. Also, that from a set of sentences $ \{ A_1, A_2, \ldots \} $ \textit{logically follows} a sentence  $ B $, means that starting from the truth of the sentences $ A_1, A_2, \ldots $ the internal semantic structure of language, not the reality the language speaks of, determines the truth of $ B $. Thus, for example, for the first-order languages the semantics of the conjunction $ \land $ determines that $ B $ logically follows from $ A \land B $. Of course, these are simple examples, and this internal language description   of logical truth and logical consequence itself is not entirely accurate. But, as with the concept of  logical symbol, in concrete and simpler situations it clearly determines whether a sentence is a logical truth, that is, whether a sentence logically follows from a set of sentences. This is, however, a good enough basis to develop a formal calculus of logical truths and logical consequence. For a first-order language without the symbol of equality and generalized quantifiers, it is easy to show by a modification of Quine's argument \cite {Qui}, Chapter 4, that this language concept of logical consequence can be described by a formal first-order logic calculus that is complete in terms of the received view. However, this proof includes Tarski’s concept of logical consequence that refers to all interpretations of the language. But a proof can be carried out without it. Namely, from the language concept of logical consequence, examining the rules of which a standard complete formal calculus is composed, it is easy to get that from the formal inference $ A_1, A_2, \ldots \vdash B $ follows  $ A_1, A_2, \ldots \models B $, in the sense of the language concept of logical consequence. The reversal can be shown by contraposition, that from $ A_1, A_2, \ldots \nvdash B $ follows $ A_1, A_2, \ldots \nvDash B $. For a standard complete formal calculus,  the condition $ A_1, A_2, \ldots \nvdash B $ is equivalent to the condition of formal consistency of the set $\{ A_1, A_2, \ldots, \lnot B \}$. Likewise, for such systems, for example for the system of classical natural deduction, Henkin-type proof of completeness (\cite{Hen}) is valid. That proof gives an interpretation which is composed of the language symbols, an interpretation that has no ontological weight, and in which the statements $ A_1, A_2 , \ldots $ are true and $ B $ is false. According to the language concept of logical consequence, this means that $A_1, A_2, \ldots \nvDash B$ is valid. This completes the proof. From this result, that in the case of first-order logic without equality, the internal language concept of logical consequence can be described by some standard formal calculus, it follows that the internal concept  coincides extensionally with the external received view concept of logical consequence. Since logical truth can be described by logical consequence, it also follows from this result that the internal language concept of logical truth can be described by some standard formal system of first-order logic without equality.

\section{Conclusion}
 Using the example of a first-order language, it is shown that understanding logic as the internal organization of language provides precise criteria for separating logical elements of language from non-logical ones. Boolean connectives and logical quantifiers are logical symbols of first-order languages, while equality and cardinal quantifiers are not.

\bibliographystyle{abbrvnat}
\bibliography{Logic}

\begin{thebibliography}{15}
\providecommand{\natexlab}[1]{#1}
\providecommand{\url}[1]{\texttt{#1}}
\expandafter\ifx\csname urlstyle\endcsname\relax
  \providecommand{\doi}[1]{doi: #1}\else
  \providecommand{\doi}{doi: \begingroup \urlstyle{rm}\Url}\fi

\bibitem[Bonnay(2008)]{Bon1}
D.~Bonnay.
\newblock {Logicality and invariance}.
\newblock \emph{Bulletin of Symbolic Logic}, 14\penalty0 (1):\penalty0 29 --
  68, 2008.

\bibitem[Bonnay(2014)]{Bon}
D.~Bonnay.
\newblock Logical constants, or how to use invariance in order to complete the
  explication of logical consequence.
\newblock \emph{Philosophy Compass}, 9\penalty0 (1):\penalty0 54--65, 2014.

\bibitem[{\v{C}}ulina(2013)]{Cul}
B.~{\v{C}}ulina.
\newblock Logic of paradoxes in classical set theories.
\newblock \emph{Synthese}, 190\penalty0 (3):\penalty0 525--547, 2013.

\bibitem[{\v{C}}ulina(2021)]{Cu0}
B.~{\v{C}}ulina.
\newblock The language essence of rational cognition with some philosophical
  consequences.
\newblock \emph{Tesis (Lima)}, 14\penalty0 (19):\penalty0 631--656, 2021.

\bibitem[Feferman(1999)]{Fef}
S.~Feferman.
\newblock {Logic, Logics, and Logicism}.
\newblock \emph{Notre Dame Journal of Formal Logic}, 40\penalty0 (1):\penalty0
  31 -- 54, 1999.

\bibitem[Gupta and Belnap(1993)]{GB}
A.~Gupta and N.~Belnap.
\newblock \emph{The Revision Theory of Truth}.
\newblock MIT, 1993.

\bibitem[Henkin(1949)]{Hen}
L.~Henkin.
\newblock The completeness of the first-order functional calculus.
\newblock \emph{Journal of Symbolic Logic}, 14\penalty0 (3):\penalty0 159--166,
  1949.

\bibitem[MacFarlane(2017)]{Mac}
J.~MacFarlane.
\newblock {Logical Constants}.
\newblock In E.~N. Zalta, editor, \emph{The {Stanford} Encyclopedia of
  Philosophy}. Metaphysics Research Lab, Stanford University, 2017.

\bibitem[Quine(1986)]{Qui}
W.~V. Quine.
\newblock \emph{Philosophy of Logic: Second Edition}.
\newblock Harvard University Press, 1986.

\bibitem[Rusnock and George(2014)]{Bol}
P.~Rusnock and R.~George.
\newblock \emph{Bernard Bolzano: Theory of Science}.
\newblock Oxford University Press, 2014.

\bibitem[Shapiro(2006)]{Sha}
S.~Shapiro.
\newblock \emph{Necessity, Meaning, and Rationality: The Notion of Logical
  Consequence}, chapter~14, pages 225--240.
\newblock John Wiley \& Sons, Ltd, 2006.

\bibitem[Sher(2012)]{Sher}
G.~Sher.
\newblock Logical quantifiers.
\newblock In D.~Graff~Fara and G.~Russell, editors, \emph{Routledge Companion
  to Philosophy of Language}, pages 579--595. Routledge, 2012.

\bibitem[Tarski(1983)]{Tar1}
A.~Tarski.
\newblock On the concept of logical consequence.
\newblock In J.~Corcoran, editor, \emph{Logic, Semantics, Metamathematics},
  pages 402--420. Hackett, 1983.

\bibitem[Tarski and Corcoran(1986)]{Tar2}
A.~Tarski and J.~Corcoran.
\newblock What are logical notions?
\newblock \emph{History and Philosophy of Logic}, 7\penalty0 (2):\penalty0
  143--154, 1986.

\bibitem[Westerståhl(2019)]{We}
D.~Westerståhl.
\newblock {Generalized Quantifiers}.
\newblock In E.~N. Zalta, editor, \emph{The {Stanford} Encyclopedia of
  Philosophy}. Metaphysics Research Lab, Stanford University, 2019.

\end{thebibliography}

\end{document}